\documentclass{aupl}

\usepackage{
amsfonts,
amsmath,
latexsym,
amssymb,
enumerate,
verbatim,
mathrsfs,
centernot,
accents,
graphicx,
rotate,
color,
epigraph,
}

\newcommand{\labbel}[1]{\label{#1} [[{\bf #1}]]}  
\renewcommand{\labbel}{\label}

\DeclareMathOperator{\diag}{Diag}

\DeclareMathOperator{\mo}{Mod}

\newtheorem{theorem}{Theorem}[section]
\newtheorem{lemma}[theorem]{Lemma}
\newtheorem{proposition}[theorem]{Proposition}

\theoremstyle{definition}
\newtheorem{definition}[theorem]{Definition}

\theoremstyle{remark}

\newtheorem{remarks}[theorem]{Remarks}

\numberwithin{equation}{section}

\begin{document}

\title{A remark about amalgamation of union of theories}

\author{Paolo Lipparini} 

\email{lipparin@axp.mat.uniroma2.it}

\urladdr{http://www.mat.uniroma2.it/~lipparin}

\address{Dipartimento di Matematica\\Viale della  Ricerca
Scientifica\\Universit\`a di Roma ``Tor Amalgamata'' \\I-00133 ROME ITALY\\ORCiD: 0000-0003-3747-6611}

\subjclass{03C52}

\keywords{Amalgamation property, existentially complete structure, compatible
theories, union of theories}


\begin{abstract}
We provide conditions under which 
the union of two first-order theories has the
amalgamation property.
\end{abstract}

\maketitle

\section{Introduction} \labbel{intro} 

The amalgamation property (AP) is an important  tool 
in algebra \cite{KMPT}, logic \cite{GM,H,MMT}
and is acquiring prominence  in computer science; see \cite{GG}
and further references there. 
In 1963 B. J{\'o}nsson \cite{J1}
 asked whether there are conditions implying AP
which can be detected just by the form of the axioms
of some theory.
While J{\'o}nsson's 
  problem seems to have  a positive solution
in a minority of cases,
there are many situations in which 
theories with AP can be merged  \cite{GG, apu} 
  in order to obtain other theories with 
AP.

In \cite{GG} some results from \cite{G} are extended in order
to show AP for union of theories.
In \cite{G,GG} a compatibility condition is assumed
with respect to some base theory which has a model-completion. 
We observe that the assumption that the base theory has a model-completion
is not necessary, actually, we do not even need to assume
that the class of models in the ``common''  
language is elementary.
The main problem still open is whether the compatibility
condition can be weakened and, if this is the case,
how far  it can be weakened. Compare \cite[p. 240]{G}.

\subsection{E.c.-compatibility} \labbel{compat}

\begin{definition} \labbel{ap}    
We assume that the reader is familiar with 
 basic notions of model theory \cite{CK,H}.

Recall that 
some class $\mathcal K$ of models  
for the same language has the \emph{amalgamation property}
(AP) if, whenever 
$\mathbf A, \mathbf B, \mathbf C \in \mathcal K$,
 $ \alpha  \colon \mathbf C \to \mathbf A$
and 
 $ \beta  \colon \mathbf C \to \mathbf B$
are embeddings, then there are a structure
$\mathbf D \in \mathcal K$ and  embeddings
$ \iota  \colon \mathbf A \to \mathbf D$
and 
 $ \eta  \colon \mathbf B \to \mathbf D$
such that 
$  \alpha  \circ \iota
=
\beta  \circ \eta$. 
\begin{equation*}
\begin{array}{lcr} 
  & \mathbf D & \cr
\enspace   \quad \quad \iota $\rotatebox[origin=c]{225}{$\dashleftarrow$}$  & & 
$\rotatebox[origin=c]{315}{$\dashleftarrow$}$  \eta   \enspace  \quad  \quad  \cr 
\enspace  \quad  \mathbf A  & &  \mathbf  B \enspace  \quad  \cr
\enspace    \quad  \quad  \alpha \nwarrow  & &    \nearrow \beta  \quad   \enspace  \quad  \cr
 & \mathbf C 
 \end{array} 
  \end{equation*}    

A quintuple $\mathbf A$, $\mathbf  B$, $\mathbf  C$, $\alpha$, $\beta$ 
as above shall sometimes be called a \emph{TBA quintuple in $\mathcal K$}, and a
model $\mathbf D$ as above shall be called an \emph{amalgamating structure}.
Reference to $\mathcal K$ shall be omitted when clear from the context.
If $\mathcal K= \mo (T)$ is the class of models of some theory $T$,
we shall say that $T$ has AP if $\mathcal K$ has AP, and similarly for
the corresponding terminology. 

If $\mathcal H$ is another class of 
models for the same language and, for every TBA quintuple in $\mathcal K$, 
there exist some $\mathbf D$ and embeddings as above, with 
$\mathbf D \in \mathcal H$, we shall say that 
$\mathcal K$ has \emph{the amalgamation property in $\mathcal H$}.  

The \emph{strong amalgamation property} means that,
under the assumption of AP,
$\mathbf D$, 
$ \iota  $
and 
 $ \eta  $ can be chosen in such a way that 
$\iota(A) \cap \eta (B) = (\alpha \circ \iota)(C)$. 
\end{definition}
 
We are going to extend some notions and results from \cite{G}
by removing the assumption that 
the ``common'' theory has model-completion, though
we need retain the amalgamation property
and work with inductive theories.
Actually, we do not need to assume that the common reduct structures
are models of some theory, we just assume  to work in some class
$\mathcal K_0$ of models for the same language $\mathscr L_0$.
Of course, the reader might always assume that 
$\mathcal K_0$ is the class of models of
some theory in $\mathscr L_0$; in this case we get the 
more usual standard notions and definitions. 

\begin{definition} \labbel{ec}  
Let $\mathcal K_0$ be a class of models in the language
$\mathscr L_0$. Recall that a model $\mathbf E$ 
is \emph{existentially complete (e.c., \emph{for short}) in $\mathcal K_0$} if 
 $\mathbf E \in \mathcal K_0$ and,

(e.c.) for every embedding $\iota:\mathbf  E \to \mathbf D$ 
into some model $\mathbf D \in \mathcal K_0$ 
and for every existential sentence $\varphi$  
in the language $\mathscr L_0 \cup E$,
if $\mathbf D \models \varphi $, then
$\mathbf  E \models \varphi $.

Here, by abuse of notation, we have identified the
model $\mathbf  E$ with the model $(\mathbf  E, e) _{e \in E} $,
that is, the model in the language $\mathscr L_0 \cup E$
obtained by adding a new constant $c_e$ to $\mathscr L_0$,
for every $e \in E$, in such a way that $c_e$ is interpreted as $e$
itself in $(\mathbf  E, e) _{e \in E} $. Similarly,
in the above condition (e.c.),  $\mathbf D$ is identified with
$(\mathbf D, \iota(e)) _{e \in E} $.
  
Of course, if $\mathcal K_0$ is closed under isomorphism,
we could equivalently suppose that $\mathbf  E$
is a submodel of  $\mathbf D$, so that there is no possible naming
conflict.   However, in what follows we shall deal 
with many embeddings at a time, hence the 
 ``substructure terminology'' will not always be viable,
mainly due to the fact that we shall not generally assume
the strong amalgamation property.
\end{definition}

The general assumptions we make are 
as follows: $T_1$  is a first order theory in the language
$\mathscr L_1$ and $\mathcal K_0$ is a class of models in the language
$\mathscr L_0$, with    
 $\mathscr L_0 \subseteq  \mathscr L_1 $.   
Subsequently, we shall also consider
a theory  $T_2$  in the language $\mathscr L_2$,
with $\mathscr L_0 = \mathscr L_1 \cap \mathscr L_2$,
We also set $\mathscr L = \mathscr L_1 \cup \mathscr L_2$
 and $T= T_1 \cup T_2$.

\begin{definition} \labbel{comp}
Suppose that $\mathcal K_0$ is a class of models
for the language $\mathscr L_0$
and $T_1$ is a first-order theory in some language
$\mathscr L_1 \supseteq \mathscr L_0$. 
 
We say that $T_1$ is
\emph{e.c.-$\mathcal K_0$-compatible} if
  \begin{enumerate}   
 \item 
For every model $\mathbf A$ of $T_1$,
the reduct $\mathbf A _{ \restriction  \mathscr L_0}   $
belongs to $\mathcal K_0$, and
\item
every model   of $T_1$
can be embedded in a model $\mathbf  B$ of $T_1$
such that  $\mathbf B _{ \restriction  \mathscr L_0}   $
is e.c.\ in $\mathcal K_0$.
  \end{enumerate}   

We could similarly give the definition of
``$\mathcal K_1$ is
e.c.-$\mathcal K_0$-compatible'', for some class 
$\mathcal K_1$ in $\mathscr L_1$, but it seems that 
the assumption is not sufficient to get significant results, namely,
here we need assume that $\mathcal K_1$   is elementary. 
 \end{definition}

\section{The main lifting lemma} \labbel{lemlem} 

\begin{lemma} \labbel{lema} 
Suppose that $\mathcal K_0$ is a class of models
for the language $\mathscr L_0$,
 $T_1$ and $T_2$   are inductive first-order theories in the languages,
respectively, $\mathscr L_1$, $\mathscr L_2$, and 
$\mathscr L_1 \cap \mathscr L_2 = \mathscr L_0$. 
Suppose further that $\mathcal K_0$ has the amalgamation property 
and both $T_1$ and $T_2$ are 
e.c.-$\mathcal K_0$-compatible.

If $\mathbf D_0 \in \mathcal K_0$, $\mathbf  E \in \mo(T_1)$,
$\mathbf  F \in \mo(T_2)$ and 
 $\iota_0 : \mathbf D_0 \to \mathbf  E  _{ \restriction \mathscr L_0 }   $,     
$\eta_0 : \mathbf D_0 \to \mathbf  F  _{ \restriction \mathscr L_0 }   $
are $\mathscr L_0$-embeddings, 
then there are a model $\mathbf  G$ of 
$T_1 \cup T_2$,  an $\mathscr L_1$-embedding
 $\iota : \mathbf  E \to \mathbf  G  _{ \restriction \mathscr L_1 }   $
and
an $\mathscr L_2$-embeddings
 $\eta : \mathbf  F \to \mathbf  G  _{ \restriction \mathscr L_2}   $
such that $\iota_0 \circ \iota = \eta _0 \circ \eta $. 
\begin{equation*}
\begin{array}{lcr} 
 & \mathbf G \cr
\mathscr L_1 $\rotatebox[origin=c]{225}{$\dashleftarrow$}$ 
 \iota & & \eta $\rotatebox[origin=c]{315}{$\dashleftarrow$}$ \mathscr L_2   \cr 
\mathbf E  & &  \mathbf  F \cr
\mathscr L_0 \nwarrow \iota_{0} & & \eta _{0} \nearrow  \mathscr L_0   \cr
 & \mathbf D_0
 \end{array}   
\end{equation*}

In addition, $ \mathbf G _{ \restriction \mathscr L_0}  $ can be taken to be  
e.c.\ in $\mathcal K_0$. 
\end{lemma} 

 \begin{proof}
Since $T_1$ is
e.c.-$\mathcal K_0$-compatible, then
there are a model $\mathbf  E_{1}$ of $T_1$
and an $\mathscr L_1$-embedding
$\iota_{1}: \mathbf  E
\to \mathbf  E_{1}$
 such that  $(\mathbf  E_{1})_{ \restriction  \mathscr L_0}   $
is e.c.\ in $\mathcal K_0$. 
Similarly,
there are a model $\mathbf  F_{1}$ of $T_2$
and an $\mathscr L_2$-embedding
$\eta _{1}: \mathbf  F
\to \mathbf  F_{1}$
 such that  $(\mathbf  F _{1})_{ \restriction  \mathscr L_0}   $
is e.c.\ in $\mathcal K_0$.
\begin{equation*}
\begin{array}{ccc} 
\mathbf E_{1}  & &  \mathbf  F_{1} \cr
\mathscr L_1 \uparrow \iota_{1} & & \eta _{1} \uparrow  \mathscr L_2   \cr
\mathbf E  & &  \mathbf  F \cr
\mathscr L_0 \nwarrow \iota_{0} & & \eta _{0} \nearrow  \mathscr L_0   \cr
 & \mathbf D_0
 \end{array}   
\end{equation*}

We now claim that 

\smallskip 

(*) there are a model $\mathbf  E_{2}$ of $T_1$,
 an $\mathscr L_1$-embedding
$\iota_{2}:   \mathbf  E_{1}
\to \mathbf  E_{2}$
and an 
$\mathscr L_0$-embedding
$\kappa _{1,2}:   \mathbf  F_{1}
\to \mathbf  E_{2}$ such that
$\iota_{0} \circ \iota_{1} \circ \iota_{2} =
 \eta _{0} \circ  \eta _{1} \circ \kappa _{1,2}$.
Furthermore, we might also assume that 
$\mathbf  E_{2}$ is e.c.\ in $\mathcal K_0$.
\begin{equation*}
\begin{array}{ccc} 
\mathbf E_{2} \cr
\mathscr L_1 \uparrow \iota_{2} &   \kappa  _{1,2 } 
$\rotatebox[origin=c]{20}{$\nwarrow$}$  \mathscr L_0  \cr
\mathbf E_{1}  & &  \mathbf  F_{1} \cr
\mathscr L_1 \uparrow \iota_{1} & & \eta _{1} \uparrow  \mathscr L_2   \cr
\mathbf E  & &  \mathbf  F \cr
\mathscr L_0 \nwarrow \iota_{0} & & \eta _{0} \nearrow  \mathscr L_0   \cr
 & \mathbf D_0
 \end{array}   
\end{equation*}

In order to prove the
first statement in (*),
 in view of the Diagram Lemma \cite[Proposition 2.1.8]{CK},
it is enough to show that  
$ \Sigma = T_1 \cup \diag _{\mathscr L_1} (\mathbf  E_{1}) 
\cup   \diag _{\mathscr L_0} (\mathbf  F_{1}) $
is a consistent theory, equivalently, by the 
compactness theorem, that every finite subset 
of $\Sigma$  
has a model.
We assume that if $d \in D_0$, then 
$(\iota_{0} \circ \iota_{1})(d)$ in  $\diag _{\mathscr L_1} (\mathbf  E_{1})$
and 
$(\eta _{0}\circ \eta _{1})(d)$ in  $\diag _{\mathscr L_0} (\mathbf  F_{1})$
correspond to the same constant in $\Sigma$,
so that the embeddings $\iota_{2}$,
$ \kappa  _{1,2}$  given by the Diagram Lemma
will commute with  $\iota_{0} \circ \iota_{1}$ and 
$\eta _{0}\circ \eta _{1}$.
(Alternatively, by taking isomorphic copies, it is no loss of 
of generality to assume that $D_0= E_1 \cap F_1$
and that  $\iota_{0} \circ \iota_{1}$ and 
$\eta _{0}\circ \eta _{1}$ are the inclusion maps.)

We shall actually show that 
$ T_1 \cup \diag _{\mathscr L_1} (\mathbf  E_{1}) 
\cup  \Gamma  $ has a model, for every 
finite subset $\Gamma$ of $\diag _{\mathscr L_0} (\mathbf  F_{1})$.
Given such a $\Gamma$, let $\gamma$ be the conjunction of the sentences in 
$\Gamma$  and let 
$ \gamma ^*=\gamma(d_1|x_1, \dots, d_k|x_k)$
be obtained from $\gamma$ by substituting the constants
$d_1, \dots $ with new variables $x_1, \dots$,  
where the $d_i$s are exactly  the constants in $\gamma$ representing 
the elements
in $  F_{1} \setminus (\eta_{0} \circ  \eta_{1})(D_0)$.
Since such constants do not appear in 
$ T_1 \cup \diag _{\mathscr L_1} (\mathbf  E_{1}) $,
then 
$ T_1 \cup \diag _{\mathscr L_1} (\mathbf  E_{1}) \cup \{  \gamma \}  $
has a model if and only if 
$ T_1 \cup \diag _{\mathscr L_1} (\mathbf  E_{1}) \cup 
\{  \exists \bar x \gamma^* \}  $
has a model.

Since $\mathcal K_0$ has AP,
there is a model $\mathbf D_0^* \in \mathcal K_0$  
which $\mathscr L_0$-amalgamates 
$\mathbf D_0$,   $ \mathbf  E_{1} $ and 
$  \mathbf  F_{1}$ relative to  $\iota_{0} \circ \iota_{1}$ 
and $\eta _{0}\circ \eta _{1}$. 
Since $\gamma$, and hence also $\gamma^*$ are quantifier-free,
then $\exists \bar x \gamma^*$ is an existential sentence.
Since $\exists \bar x \gamma^*$
is an $\mathscr L_0 \cup D_0$-sentence
which holds in  $  \mathbf  F_{1}$
(recall that elements of $D_0$ and of 
$(\eta _{0}\circ \eta _{1})(D_0)$ are given the same name),
then $\exists \bar x \gamma^*$ holds in $\mathbf D_0^*$, as well.
Since $(\mathbf  E_{1})_{ \restriction  \mathscr L_0}   $
is e.c.\ in $\mathcal K_0$ and 
$\mathscr L_0 \cup D_0 $-embeds in $\mathbf D_0^*$,
then  $\exists \bar x \gamma^*$ holds in 
$\mathbf  E_{1}$. Thus 
$\mathbf  E_{1}$ is a model of
$ T_1 \cup \diag _{\mathscr L_1} (\mathbf  E_{1}) \cup 
\{  \exists \bar x \gamma^* \}  $.

As mentioned above, applying the above argument to each
finite subset $\Gamma$ of $\diag _{\mathscr L_0} (\mathbf  F_{1})$
we show that $ \Sigma = T_1 \cup \diag _{\mathscr L_1} (\mathbf  E_{1}) 
\cup   \diag _{\mathscr L_0} (\mathbf  F_{1}) $
has a model, and this implies the first statement in (*).
The second statement in (*) is then immediate 
from the assumption that $T_1$ is 
e.c.-$\mathcal K_0$-compatible: just use compatibility
in order to embed the model constructed in the proof of the first statement
into some model e.c.\ in $\mathcal K_0$.

Having proved (*), then, by a similar and slightly 
simpler argument, we get that there are 
a model $\mathbf  F_{3}$ of $T_2$,
 an $\mathscr L_2$-embedding
$\eta _{3}:   \mathbf  F_{1}
\to \mathbf  F_{3}$
and an 
$\mathscr L_0$-embedding
$\kappa _{2,3}:   \mathbf  E_{2}
\to \mathbf  F_{3}$ such that
$\eta _{3}  = \kappa _{1,2} \circ \kappa _{2,3}$.
\begin{equation*}
\begin{array}{ccc} 
  & &  \mathbf  F_{3} \cr
 &   \kappa  _{2,3 } 
$\rotatebox[origin=c]{-20}{$\nearrow$}$
 \mathscr L_0  \cr
\mathbf E_{2} & & \eta_{3} \uparrow \mathscr L_2  \cr
\mathscr L_1 \uparrow \iota_{2} &   \kappa  _{1,2 } 
$\rotatebox[origin=c]{20}{$\nwarrow$}$
 \mathscr L_0  \cr
\mathbf E_{1}  & &  \mathbf  F_{1} \cr
\mathscr L_1 \uparrow \iota_{1} & & \eta _{1} \uparrow  \mathscr L_2   \cr
\mathbf E  & &  \mathbf  F \cr
\mathscr L_0 \nwarrow \iota_{0} & & \eta _{0} \nearrow  \mathscr L_0   \cr
 & \mathbf D_0
 \end{array}   
\end{equation*}

Indeed, in view of the Diagram Lemma,
it is enough to show that  every finite subset of 
$ T_2 \cup \diag _{\mathscr L_2} (\mathbf  F_{1}) 
\cup   \diag _{\mathscr L_0} (\mathbf  E_{2}) $
has a model.
Arguing as above, it is enough to consider,
in place of $\diag _{\mathscr L_0} (\mathbf  E_{2})$, some
existential sentence $\exists \bar x \gamma^*$
  in $\mathscr L_0 \cup   F_1$. 
Any such sentence holds in $\mathbf  E_{2}$, hence it holds in
$\mathbf  F_{1}$
since  $\mathbf  F_{1}$ is 
is e.c.\ in $\mathcal K_0$.
Hence $\mathbf  F_{1}$ is a model 
of 
$ T_2 \cup \diag _{\mathscr L_2} (\mathbf  F_{1}) 
\cup  \{  \exists \bar x \gamma^*\} $.
The existence of $\mathbf  F_{3}$ together with appropriate
embeddings follows.

Having proved the existence of some $\mathbf  F_{3}$,
 then, by taking an appropriate extension,
we can further assume that $\mathbf  F_{3}$ is e.c.\ in $\mathcal K_0$,
 since $T_2$ is 
e.c.-$\mathcal K_0$-compatible.

Iterating the above arguments, we get chains of models
with a commutative diagram of embeddings
\begin{equation*}
\begin{array}{ccc}
\vdots && \vdots \cr
\mathbf E_{6}  & & \eta_{7} \uparrow \mathscr L_2 \cr
 &   \kappa  _{5,6 } 
$\rotatebox[origin=c]{20}{$\nwarrow$}$
  \mathscr L_0  \cr
\mathscr L_1 \uparrow \iota_{6} & &  \mathbf  F_{5} \cr 
 &   \kappa  _{4,5 } $\rotatebox[origin=c]{-20}{$\nearrow$}$
 \mathscr L_0  \cr
\mathbf E_{4}  & & \eta_{5} \uparrow \mathscr L_2 \cr
 &   \kappa  _{3,4 } $\rotatebox[origin=c]{20}{$\nwarrow$}$  \mathscr L_0  \cr
\mathscr L_1 \uparrow \iota_{4} & &  \mathbf  F_{3} \cr
 &   \kappa  _{2,3 } $\rotatebox[origin=c]{-20}{$\nearrow$}$
  \mathscr L_0  \cr
\mathbf E_{2} & & \eta_{3} \uparrow \mathscr L_2  \cr
\mathscr L_1 \uparrow \iota_{2} &   \kappa  _{1,2 } 
$\rotatebox[origin=c]{20}{$\nwarrow$}$  \mathscr L_0  \cr
\mathbf E_{1}  & &  \mathbf  F_{1} \cr
\mathscr L_1 \uparrow \iota_{1} & & \eta _{1} \uparrow  \mathscr L_2   \cr
\mathbf E  & &  \mathbf  F \cr
\mathscr L_0 \nwarrow \iota_{0} & & \eta _{0} \nearrow  \mathscr L_0   \cr
 & \mathbf D_0
 \end{array}   
\end{equation*}
Since $T_1$ is inductive, the direct limit $\mathbf E _{\infty } $ 
of the models $\mathbf E$, 
$\mathbf E_{1}$, $\mathbf E_{2i}$ through the embeddings 
$\iota_1$, $\iota_{2i}$ is a model of $T_1$.
Similarly, the direct limit $\mathbf F _{\infty } $ 
of the models $\mathbf F$, $\mathbf  F_{2i+1}$ through the embeddings 
$\eta_{2i+1}$ is a model of $T_2$.
The embeddings $\kappa _{i, i+1} $ 
  induce an $\mathscr L_0$-isomorphism $\kappa$ between 
$\mathbf  F_{\infty } $ and $\mathbf  E_{\infty } $, so that if  
we use $\kappa$ to transport the 
$\mathscr L_2 \setminus \mathscr L_0$
structure of $\mathbf  F_{\infty } $ into $\mathbf E_{\infty } $, we get a model
$\mathbf G$
of $T_1 \cup T_2$. This is possible since we have assumed that
  $\mathscr L_1 \cap \mathscr L_2 = \mathscr L_0$. 

The model $ \mathbf G _{ \restriction \mathscr L_0}  $ is  
e.c.\ in $\mathcal K_0$, since it is a limit of 
models which are e.c.\ in $\mathcal K_0$.
The embeddings $\iota_1$ and $\iota_{2i}$
provide an embedding  $\iota$ of 
$\mathbf  E   $  
into $ \mathbf  G _{ \restriction  \mathscr L_1} = \mathbf E_{\infty }$.
The embeddings  $\eta_{2i+1}$ and $\kappa$ 
provide an embedding $\eta$  of 
$\mathbf  F  $  
into  $ \mathbf  G _{ \restriction  \mathscr L_2}$.
The embeddings $\iota$ and $\eta$ agree on $\mathscr L_0 \cup D_0$,
since the above diagrams are commutative, hence we get 
the conclusion.
 \end{proof}

\section{Compatibility and AP for union of theories} \labbel{mod} 

For universal theories, the following theorem
 generalizes \cite[Proposition 4,2]{G},
to the effect that, as we mentioned, we need not assume that 
 $\mathcal K_0$ is an elementary class with model-completion.
We still assume that $\mathcal K_0$ has the amalgamation property,
a consequence of the existence of the model-completion.

\begin{theorem} \labbel{compu} 
Suppose that $\mathcal K_0$ is a class of models
for the language $\mathscr L_0$,
 $T_1$ and $T_2$   are inductive first-order theories in the languages,
respectively, $\mathscr L_1$, $\mathscr L_2$, and 
$\mathscr L_1 \cap \mathscr L_2 = \mathscr L_0$. 
 
If $\mathcal K_0$ has the amalgamation property 
and both $T_1$ and $T_2$ are 
e.c.-$\mathcal K_0$-compatible,
then
$T_1 \cup T_2$ is
e.c.-$\mathcal K_0$-compatible.
 \end{theorem}

\begin{proof}
Let $\mathbf  C$ be a model of $T_1 \cup T_2$.
Apply Lemma \ref{lema} taking 
$\mathbf D_0 = \mathbf  C_{ \restriction  \mathscr L_0}$,
$\mathbf  E= \mathbf  C_{ \restriction  \mathscr L_1}$
and $\mathbf  F = \mathbf  C_{ \restriction  \mathscr L_2}$,
with $\iota_0$ and $\eta_0$ the identity function.
Since $\iota_0$ and $\eta_0$ are the identity,
then the functions $\iota$ and $\eta$ given by Lemma \ref{lema}
coincide, so they furnish  an $\mathscr L_1 \cup \mathscr L_2$-embedding
of $\mathbf  C$ 
into $\mathbf  G$, whose $\mathscr L_0$ reduct 
is e.c.\ in $\mathcal K_0$.
\end{proof}
 
Notice that, under the assumptions in Theorem \ref{compu},
 it might happen that $T_1 \cup T_2$
has no model. Formally, the theorem remains 
true in this trivial case, due to the usual conventions
about universal quantification over empty classes
(if $T$ is contradictory, then $T$ is vacuously  
e.c.-$\mathcal K_0$-compatible).

\begin{definition} \labbel{common}
Suppose that  $T_1$ and $T_2$   are  first-order theories in the languages,
respectively, $\mathscr L_1$, $\mathscr L_2$, 
set 
$\mathscr L= \mathscr L_1 \cup \mathscr L_2 $,
$\mathscr L_0= \mathscr L_1 \cap \mathscr L_2 $
and suppose that $\mathcal K_0$ is a class of models
for the language $\mathscr L_0$.

We say that  $T_1$ and $T_2$  
have \emph{common subcompatible amalgamation 
over $\mathcal K_0$} if,  
 for every TBA-quintuple $\mathbf A$, $\mathbf  B$, $\mathbf  C$,
$\alpha$, $\beta$ in 
$T_1 \cup T_2$,
there are a model $\mathbf D_0 \in \mathcal K_0$,
models $\mathbf  E \in \mo ( T_1)$ and  $\mathbf  F \in \mo ( T_2)$,
$\mathscr L_0$-embeddings 
$\alpha_1: \mathbf  A \to \mathbf D_0$,
$ \beta_1 :  \mathbf  B \to \mathbf D_0$,
$\iota_0: \mathbf D_0 \to \mathbf  E$ and
$\eta_0: \mathbf D_0 \to \mathbf  F$
such that $\alpha \circ \alpha _1 = \beta \circ \beta_1 $
and moreover 
$ \alpha _1 \circ \iota_0$ and
$\beta_1 \circ \iota_0$ are $\mathscr L_1$-embeddings
and
$ \alpha _1 \circ \eta _0$ and
$\beta_1 \circ \eta _0$ are $\mathscr L_2$-embeddings.    
\begin{gather*}
\begin{array}{lcr} 
\enspace \quad \mathbf E  & &  \mathbf  F \enspace \quad \cr
\enspace  \quad  \quad  \mathscr L_0 \nwarrow \iota_0 & & \eta_0 \nearrow  \mathscr L_0 \enspace  \quad  \quad   \cr 
\mathscr L_1 \uparrow  & \mathbf D_0 & \uparrow \mathscr L_2\cr
\enspace  \quad  \quad  \mathscr L_0 \nearrow \alpha _1 & & 
\beta_1 \nwarrow  \mathscr L_0   \enspace  \quad  \quad  \cr 
\enspace  \quad  \mathbf A  & &  \mathbf  B \enspace  \quad  \cr
\enspace  \quad  \quad  \mathscr L \nwarrow \alpha  & & \beta   \nearrow  \mathscr L   \enspace  \quad  \quad  \cr
 & \mathbf C 
 \end{array} 
\\
\text{with $\mathscr L= \mathscr L_1 \cup \mathscr L_2$, 
$\mathscr L_0= \mathscr L_1 \cap \mathscr L_2$,}
\\
\text{furthermore, 
$\beta_1 \circ \iota_0$ an $\mathscr L_1$-embedding,
$ \alpha _1 \circ \eta _0$ an
 $\mathscr L_2$-embedding.}   
\end{gather*}
 \end{definition}   

The above notion is slightly more general than the notions introduced
in \cite[Section 3.1]{GG} (compare the displayed formula
in Theorem  3.1 therein). Here we are not assuming that theories
are universal. 

\begin{remarks} \labbel{rms}   
(a) There is a case in which the above notion can be simplified.
Suppose that $\mathcal K_0$ has pushouts, hence AP implies that
a pushout is an amalgamating structure.  
Let us say that $T_1$ has \emph{amalgamation over 
$\mathcal K_0$-pushouts} if,
 for every TBA quintuple $\mathbf A$, $\mathbf  B$, $\mathbf  C$,
$\alpha$, $\beta$   
of models of $T_1$, there is an amalgamating model 
whose $\mathscr L_0$ reduct  extends the pushout in $\mathcal K_0$
of $\mathbf A _{ \restriction \mathscr L_0}  $ and
 $\mathbf  B_{ \restriction \mathscr L_0}  $ over
 $\mathbf  C_{ \restriction \mathscr L_0} $.     

With the above definition, 
if both $T_1$ and $T_2$ have amalgamation over 
$\mathcal K_0$-pushouts,
then trivially $T_1$ and $T_2$  
have common subcompatible amalgamation 
over $\mathcal K_0$, since we can take the above pushout 
as $\mathbf  D_0$ in Definition \ref{common}.  

The notion of  amalgamation over 
$\mathcal K_0$-pushouts has implicitly appeared in many
situations.

(a1) if $\mathcal K_0$ is the class of structures in the 
empty language, and, for simplicity, modulo isomorphism,
 we suppose that
$\mathbf  C \subseteq \mathbf A$,  
$\mathbf  C \subseteq \mathbf  B$
and $C =A \cap B $, 
then the pushout 
of the above structures
is $A \cup B$ (with no structure).  
 
In this case, $T_1$ has amalgamation over 
$\mathcal K_0$-pushouts
exactly if $T_1$ has the strong amalgamation property.

(a2) If $\mathcal K_0$ is the class of partially
ordered sets (posets, for short), then  
amalgamation over 
$\mathcal K_0$-pushouts
is equivalent to  a notion called the 
\emph{superamalgamation property}.
The standard proof of AP for  posets
is obtained  by constructing the pushout
(again, over $A \cup B$),
then
amalgamation over 
$\mathcal K_0$-pushouts
translates immediately into the superamalgamation property.
See \cite{apu} for more about the above notions.

(b) There are cases in which  common subcompatible amalgamation
alone implies that $T_1 \cup T_2$ has AP. See Proposition \ref{oss}.
In the following theorem we shall show that
$T_1 \cup T_2$ has AP when common subcompatible amalgamation
holds and both $T_1$ and $T_2$ are 
e.c.-$\mathcal K_0$-compatible.

(c)
Clearly, if $T_1 \cup T_2$ has AP, then $T_1$ and $T_2$  
have common subcompatible amalgamation 
over $\mathcal K_0$, since 
if $T_1 \cup T_2$ has AP we can take 
$E=D=F$ in the above definition 
(this observation implies that the assumption
of common subcompatible amalgamation is necessary in the next
theorem.)
On the other hand, formally, 
 common subcompatible amalgamation  for 
$T_1$ and $T_2$ does not imply AP for 
$T_1$ or $T_2$:
 common subcompatible  amalgamation
only implies that, for example,  the class $\mathcal K_1$
has AP into $\mo(T_1)$, where  $\mathcal K_1$ is
the class 
 of those models of $T_1$ 
which can be expanded to some model 
of $T_2$ (indeed, $\mathbf E $ is in  $\mo(T_1)$
and amalgamates $\mathbf A$, $\mathbf  B$ and $\mathbf  C$,
which are in $\mathcal K_1$. This is due to the request that  
$\iota \circ \iota_0$ and
$\eta \circ \iota_0$ are $\mathscr L_1$-embeddings). 
\end{remarks} 

\begin{theorem} \labbel{ama} 
Suppose that $\mathcal K_0$ is a class of models
for the language $\mathscr L_0$,
 $T_1$ and $T_2$   are inductive first-order theories in the languages,
respectively, $\mathscr L_1$, $\mathscr L_2$, and 
$\mathscr L_1 \cap \mathscr L_2 = \mathscr L_0$. 
 Suppose further that $\mathcal K_0$,
 has the amalgamation property, 
that both $T_1$ and $T_2$ are 
e.c.-$\mathcal K_0$-compatible and 
that $T_1$ and $T_2$  
have common subcompatible amalgamation 
over $\mathcal K_0$.

Then $T_1 \cup T_2$ has the amalgamation property.
 \end{theorem}

\begin{proof}
Given a TBA quintuple $\mathbf A$, $\mathbf  B$, $\mathbf  C$,
$\alpha$, $\beta$ in $T_1 \cup T_2$, 
the assumption of common subcompatible amalgamation
provides models 
 $\mathbf D_0 \in \mathcal K_0$,
 $\mathbf  E \in \mo(T_1)$, $\mathbf  F \in \mo(T_2)$
 and embeddings
 $\iota_0   $     
and $\eta_0   $
satisfying the assumptions of Lemma \ref{lema}. 
 
Then apply Lemma \ref{lema}
getting
 a model $\mathbf  G$ of 
$T_1 \cup T_2$,  an $\mathscr L_1$-embedding
 $\iota : \mathbf  E \to \mathbf  G  _{ \restriction \mathscr L_1 }   $
and
an $\mathscr L_2$-embeddings
 $\eta : \mathbf  F \to \mathbf  G  _{ \restriction \mathscr L_2}   $
such that $\iota_0 \circ \iota = \eta _0 \circ \eta $. 

Since $\alpha_1 \circ \iota_0$ is an $\mathscr L_1$ embedding,
then   $ \alpha_1 \circ \iota_0 \circ \iota$
 is an $\mathscr L_1$ embedding.
Since $\alpha_1 \circ \eta_0$ is an $\mathscr L_2$ embedding,
then   $\alpha_1 \circ \eta_0 \circ \eta$
 is an $\mathscr L_2$ embedding.
Since $\iota_0 \circ \iota = \eta _0 \circ \eta $,
then  $ \alpha_1 \circ \iota_0 \circ \iota= \alpha_1 \circ \eta_0 \circ \eta$,
henceforth
$\alpha ^*=
 \alpha_1 \circ \iota_0 \circ \iota$  is an $\mathscr L_1 \cup \mathscr L_2$ embedding.

Similarly, $ \beta ^*=
 \beta _1 \circ \iota_0 \circ \iota= \beta _1 \circ \eta_0 \circ \eta$
is an $\mathscr L_1 \cup \mathscr L_2$ embedding.
Finally, $ \alpha \circ \alpha ^*=
 \alpha \circ \alpha_1 \circ \iota_0 \circ \iota=
 \beta  \circ \beta _1 \circ \iota_0 \circ \iota=
\beta \circ \beta ^*$, hence $\mathbf G$ amalgamates 
$\mathbf A$ and $\mathbf  B$ over $\mathbf  C$ 
by means of $\alpha^*$ and $\beta^*$. 
\end{proof}

\section{Further remarks} \labbel{fur}

 \begin{proposition}  \labbel{oss}
As usual by now, 
suppose that  $T_1$ and $T_2$   are  theories in the languages
 $\mathscr L_1$, $\mathscr L_2$, with
$\mathscr L_0= \mathscr L_1 \cap \mathscr L_2 $
and  $\mathcal K_0$ is a class of models
for the language $\mathscr L_0$.
Suppose further that both 
$T_1$ and $T_2$   
are universal theories.

The assumption that 
$T_1$ and $T_2$  
have common subcompatible amalgamation 
over $\mathcal K_0$
is a necessary and sufficient condition
for 
$T_1 \cup T_2$ to have AP 
in each of the following cases.
  \begin{enumerate}[(a)]    
\item 
The languages $\mathscr L_0$, $\mathscr L_1$ and $\mathscr L_2$
  contain no function symbol of arity $\geq 2$.
\item
Both $\mathscr L_1 \setminus \mathscr L_0$ and
$\mathscr L_2 \setminus \mathscr L_0$ are  relational.

\item
In the definition of 
common subcompatible amalgamation we can take
 $\mathbf D_0$  to be the push-out in $\mathscr L_0$
of $\mathbf A _{ \restriction \mathscr L_0 }  $ and 
$\mathbf  B_{ \restriction \mathscr L_0 }$
 over $\mathbf  C_{ \restriction \mathscr L_0 }$,
moreover, 
 $\mathscr L_1 \setminus \mathscr L_0$ and 
$\mathscr L_2 \setminus \mathscr L_0$
contain no function symbol of arity $\geq 2$ and, finally,
$T_1$ and $T_2$ assert that every
function in  $\mathscr L_1 \setminus \mathscr L_0$,
respectively,
$\mathscr L_2 \setminus \mathscr L_0$
is an $\mathscr L_0$ endomorphism.
 \end{enumerate} 
 \end{proposition}   

\begin{proof}
As noticed in Remark \ref{rms}(c), 
necessity holds with no further assumption. 

(a) Let $\mathbf D_0$, $\mathbf  E$, $\mathbf  F$, 
$\iota_0$ and $\eta_0$ be given by  common subcompatible amalgamation,
Definition \ref{common}.
\begin{gather*}
\begin{array}{lcr} 
\enspace \quad \mathbf E  & &  \mathbf  F \enspace \quad \cr
\enspace  \quad  \quad  \mathscr L_0 \nwarrow \iota_0 & & \eta_0 \nearrow  \mathscr L_0 \enspace  \quad  \quad   \cr 
\mathscr L_1 \uparrow  & \mathbf D_0 & \uparrow \mathscr L_2\cr
\enspace  \quad  \quad  \mathscr L_0 \nearrow \alpha _1 & & 
\beta_1 \nwarrow  \mathscr L_0   \enspace  \quad  \quad  \cr 
\enspace  \quad  \mathbf A  & &  \mathbf  B \enspace  \quad  \cr
\enspace  \quad  \quad  \mathscr L \nwarrow \alpha  & & \beta   \nearrow  \mathscr L   \enspace  \quad  \quad  \cr
 & \mathbf C 
 \end{array} 
\\
\text{with $\mathscr L= \mathscr L_1 \cup \mathscr L_2$, 
$\mathscr L_0= \mathscr L_1 \cap \mathscr L_2$,}
\\
\text{$\beta_1 \circ \iota_0$ an $\mathscr L_1$-embedding,
$ \alpha _1 \circ \eta _0$ an
 $\mathscr L_2$-embedding.}   
\end{gather*}

Set $D^* = \alpha _1 (A) \cup \beta  _1 (B)$.
Since no function of arity $\geq 2$
is present, then $D^* $
can be turned to an $\mathscr L_1 \cup \mathscr L_2$  model
by setting $f( \alpha _1 (a))= \alpha _1 (f _{ \mathbf A} (a))$
for $a \in A$ and 
$f( \beta _1 (b))= \beta  _1(f_{ \mathbf B}(b))$
for $b \in B$, for every $f$ in the language
$\mathscr L_1 \cup \mathscr L_2$.
The interpretation of relation symbols
is induced by $\iota_0$ and $ \eta _0$, namely, if, say,
$R \in \mathscr L_1$, let $R(d_1, d_2, \dots)$
hold in $D^*$ if and only if   
$R(\iota_0(d_1), \iota_0(d_2), \dots)$
holds in $\mathbf E$.

We need to check that the above definition is correct.
First, if $a \in A$, then $f_{ \mathbf A}(a) \in A$,  
hence $ \alpha _1(f_{ \mathbf A}(a))$ does belong to 
$\alpha _1 (A) $, which is a subset of $D^* $.
If $\alpha _1 (a)= \beta  _1 (b)$,
for $a \in A$ and $b \in B$, and, say, $f \in \mathscr L_1$, then  
$ (\iota_0 \circ \alpha _1) (a)= (\iota_0 \circ  \beta  _1 )(b)$,
hence
$\iota_0 ( \alpha _1 (f(a))) =
 f((\iota_0 \circ \alpha _1) (a))=f( (\iota_0 \circ  \beta  _1 )(b))=
\iota_0 ( \beta  _1 (f(b)))$, since 
$ \iota_0 \circ \alpha _1 $ and $ \iota_0 \circ  \beta  _1 $
are $\mathscr L_1$ homomorphisms. Since 
$\iota_0$ is injective, then  $\alpha _1 (f(a))=\beta  _1 (f(b))$. 
This means that $f$ is well-defined.
The argument also shows that 
$\iota_0$ is an embedding from $\mathbf D^* _{ \restriction \mathscr L_1}  $
into  $\mathbf E$, thus $\mathbf D^*$ is a model of 
$T_1$, since $T_1$ is universal. Similarly,     
 $\mathbf D^*$ is a model of 
$T_2$.

By the very definition of 
functions in $\mathbf D^*$,
we have that $\alpha_1$ is 
an $\mathscr L_1 \cup \mathscr L_2$
embedding, as far as functions are concerned.
As for relations, if, say, $R \in \mathscr L_1$, then 
$R _{ \mathbf A} (a_1, a_2, \dots )$   
if and only if $R _{ \mathbf E} (( \alpha _1 \circ \iota_0)(a_1),
 ( \alpha _1 \circ \iota_0)(a_2), \dots )$, since 
$\alpha _1 \circ \iota_0$ is an $\mathscr L_1$ embedding,
if and only if $R _{ \mathbf D^*} (\alpha _1 (a_1),
 \alpha _1 (a_2), \dots )$, by definition.
Similarly, 
$\alpha_1$ is 
an $\mathscr L_1 \cup \mathscr L_2$
embedding, hence $\mathbf D^*$
amalgamates $\mathbf A$ and $\mathbf  B$
over $\mathbf  C$.  

(b) is much easier,  $\mathbf D_0$ can be expanded
to a substructure of both  $\mathbf  E$
(in $\mathscr L_1$)  and $\mathbf  F$ (in $\mathscr L_2$),
thus becoming a model of $T_1 \cup T_2$, since both 
$T_1  $ and $ T_2$ are universal.

(c) If $f \in \mathscr L_1 \setminus \mathscr L_0$,
then $f$ can be thought of as an $\mathscr L_0$ endomomorphism,
for every  model of $T_1$; in particular, this applies
to $\mathbf A _{ \restriction \mathscr L_0}  $,
 $\mathbf  B_{ \restriction \mathscr L_0}$ and  
$\mathbf  C_{ \restriction \mathscr L_0}$.
Since   $\mathbf A _{ \restriction \mathscr L_0}  $ and
 $\mathbf  B_{ \restriction \mathscr L_0}$   
 embed into 
$\mathbf D_0$ over $\mathbf  C_{ \restriction \mathscr L_0}$,
then $f_{ \mathbf A}$ can be also thought of as an
homomorphism from $\mathbf A$ to $\mathbf D_0$, and the same
for $f_{ \mathbf B}$. Since $\mathbf D_0$ is a pushout,
$f_{ \mathbf A}$ and $f_{ \mathbf B}$ extend  uniquely 
 to a 
homomorphism from $\mathbf D_0$ to $\mathbf D_0$.
This means that $f$ can be extended uniquely on $\mathbf D_0$,
under the assumption that the diagram remains commutative. 
 See \cite[Section 5]{apuvar} for further details. 
The facts that in this way $\alpha_1$ becomes an $\mathscr L_1$ embedding   
and that $\alpha_1 \circ \iota_0$ is assumed to be  an $\mathscr L_1$ embedding 
imply that $\iota_0$ is an  $\mathscr L_1$ embedding,
arguing as in (a) above.  

Symmetrically, every  $f \in \mathscr L_2 \setminus \mathscr L_0$
can be extended to a function in $\mathbf D_0$, so that $\mathbf D_0$
becomes a model of $T_1 \cup T_2$, since $\mathbf D_0$ becomes
an $\mathscr L_1$ substructure of $\mathbf  E$ and $T_1$
is universal, and similarly for $\mathscr L_2$ and $T_2$.  
\end{proof}
  
See \cite{apu} for other results and arguments similar 
to Proposition \ref{oss}. 
 
\smallskip

This is a preliminary version, it might contain
inacuraccies.
We have not yet performed a completely accurate search in order to
check whether some of the results presented here are already known.
Credits for already known results should go to the original discoverers.

\end{document}